# A note on $m_h(A_k)$

**History**

| | | |
|---|---|---|
| 0.00 | 20-01-00 | Started: most of the theoretical background completed. |
| 0.01 | 24-01-00 | Everything complete except notes about the algorithms used. |
| 1.00 | 01-02-00 | First version complete |
| 1.01 | 21-09-14 | Add Abstract and References prior to ArxIv publication; no other changes |

**Further work**

Look for a set which has more than one d(h) > d(h-1) for different values of h.
Look for a set with the simplest possible "pattern" for making d(h) > d(h-1).
Investigate why, overall, there seem to be so few h-gaps (perhaps because most sets have $h_2 <= h_0$).

**Abstract**

$A_k = \{1, a_2, ... a_k\}$ is an h-basis for n if every positive integer not exceeding n can be expressed as the sum of no more than h values $a_i$; we write $n = n_h(A_k)$. An extremal h-basis $A_k$ is one for which n is as large as possible, and then we write $n = n_h(k)$. Computing such extremal bases has become known as the "global" Postage Stamp Problem.

The "local" Postage Stamp Problem is concerned with properties of particular sets $A_k$, and it is clear that sets where $n_h(A_k)$ does not exceed $a_k$ are of little interest. We define $h_0(k)$ to be the smallest value of h for which $n_h(A_k)$ exceeds $a_k$; such sets are called "admissible".

We say that a value n can be "generated" by $A_k$ if it can be expressed as no more than h values $a_i$, or - equivalently - if it can be expressed as the sum of exactly h values $a_i$ from the set $A'_k = \{0, 1, a_2 ... a_k\}$. No values greater than $ha_k$ can be generated, and we now consider the number of values less than $ha_k$ that have no generation, denoted $m_h(A_k)$ - essentially a count of the number of "gaps" (see Challis [1], and Selmer [5] page 3.1).

It is easy to show that for some value $h_2(k)$ exceeding $h_0(k)$ the difference $m_h(A_k) - m_{h+1}(A_k)$ remains constant - that is, the "pattern" of missing values between $ha_k$ and $(h+1)a_k$ does not change as h increases. Here we are interested in the pattern of missing values for values that lie between $h_0$ and $h_2$.

On page 7.8 of Selmer [5] he conjectures that the sequence of differences $m_h(A_k) - m_{h+1}(A_k)$ is non-increasing as h runs from $h_0$ to $h_2$. When I came across this conjecture I could not convince myself that it was likely to be true, having found a possible error in Selmer's justification. I wrote to him in November 1995, and early in 1996 he replied, agreeing that this might be the case and hoping that I might be able to find a counter-example. This paper records my successful search for a counter example, eventually found late in 1999.

**Summary**

This document describes experiments to investigate Selmer's conjecture that the number of gaps filled in as h increases from $h_0 - 1$ to $h_2$ decreases monotically; the conjecture is shown to be false by counter-example.

**Notation**

For a given set $A_k$:

$n(h, A_k)$ is the set's h-range, or cover.

$h_0$ is the smallest value of h such that $n(h, A_k) >= a_k$.

$h_1$ is the smallest value of h such that $n(h+1, A_k) = n(h, A_k) + a_k$ for all $h >= h_1$.

$n(h, A_k) < x < ha_k$ is a *gap at level h* if x has no h-representation.

$h_2$ is the smallest value of h such that for all $h \geq h_2$ if x is a gap at level h then $x + a_k$ is a gap at level h+1.

$m(h, A_k)$ is the number of gaps at level h.

$d(h, A_k) = m(h-1, A_k) - m(h, A_k)$     [Note: Selmer defines this as $m(h, A_k) - m(h+1, A_k)$]

$n(h_0-1, A_k) < x < (h_0-1)a_k$ is an *m-gap* if $x_h$ has an h-representation for h = m, but no h-representation for h < m, where $x_h = x + (h - (h_0-1))a_k$.

Since we are usually talking about a single set $A_k$, we normally write n(h), m(h), d(h) in place of $n(h, A_k)$ etc.

Selmer's conjecture is that:

$$d(h_0) \geq d(h_0+1) \geq \ldots \geq d(h_2) \geq d(h_2+1) = 0$$

We show by counter-example that this is not so.

## $h_0$, $h_1$ and $h_2$

*Theorem 1*

    $n(h+1) \geq n(h) + a_k$   for all   $h \geq h_0-1$

*Proof*

    Every value $0 \leq x \leq n(h)$ has an h-representation, so every value $a_k \leq x \leq n(h) + a_k$ has an (h+1)-representation: simply use one more addend $a_k$.

    Every value $0 \leq x \leq a_k$ has an $h_0$-representation (and hence an (h+1)-representation for $h \geq h_0-1$) by definition, thus proving the theorem.

*Theorem 2*

    There exists a value $h_1 \geq h_0-1$ such that $n(h+1) = n(h) + a_k$ for all $h \geq h_1$.

*Proof*

    Consider $D(h) = n(h) - ((h+1)a_{k-1} - a_k)$.

    For each increment in h, n(h) increases by at least $a_k$ (by Theorem 1) and $((h+1)a_{k-1} - a_k)$ increases by exactly $a_{k-1}$; so D(h) increases as h increases, and there must be a minimal value $h_s \geq h_0-1$ such that $D(h_s) \geq 0$:

    ie     $n(h) \geq (h+1)a_{k-1} - a_k$      for all   $h \geq h_s$                  (1)

We now show that $n(h) + a_k + 1$ has no (h+1)-representation for $h \geq h_s$; but Theorem 1 says that $n(h+1) \geq n(h) + a_k$, and so $n(h+1) = n(h) + a_k$ for all $h \geq h_s$ and the theorem is proved.

Let $x = n(h) + a_k + 1$ for $h \geq h_s$, and suppose it has the (h+1)-representation:

    $x = c_k a_k + c_{k-1} a_{k-1} + \ldots + c_1 a_1$

If $c_k > 0$, then $x' = x - a_k = n(h) + 1$ has the h-representation:

    $x' = (c_k-1)a_k + c_{k-1}a_{k-1} + \ldots + c_1 a_1$

which is contrary to our assumption that n(h) is the h-range.

So the (h+1)-representation of x must have the form:

$$x = c_{k-1}a_{k-1} + ... + c_1a_1$$

and so $x <= (h+1)a_{k-1}$; so:

$$x = n(h) + a_k + 1 <= (h+1)a_{k-1} \implies n(h) < (h+1)a_{k-1} - a_k$$

which contradicts our assumption that $h >= h_s$ by (1).

*Theorem 3*

There exists $h_2 >= h_0 - 1$ such that for all $h >= h_2$ if x is a gap at level h, then $x + a_k$ is a gap at level h+1.

*Proof*

Let $h >= h_s$ and let $n(h) < x < ha_k$ be a gap at level h: that is, a value with no h-representation.

Suppose $x' = x + a_k$ is not a gap at level h+1: that is, x' has an (h+1)-representation:

$$x' = c_k a_k + ... + c_1 a_1$$

If $c_k > 0$, then $x = x' - a_k$ has an h-representation, and so is not a gap at level h; so $c_k$ must be zero and hence $x' <= (h+1)a_{k-1}$.

Now $n(h) < x = x' - a_k <= (h+1)a_{k-1} - a_k$, which contradicts our assumption that $h >= h_s$ by (1).

*Theorem 4*

$h_2 >= h_1$

*Proof*

If x is the first gap at level h, then the h-range $n(h) = x - 1$; so Theorem 2 can be equivalently stated as:

For all $h >= h_1$, if x is the smallest gap at level h, then $x + a_k$ is the smallest gap at level h+1.

In other words, for $h >= h_1$ the *first* gap persists, whereas for $h >= h_2$ *all* gaps persist; hence $h_2 >= h_1$.

We say that the h-range is *stabilised* for $h >= h_1$, and that *pattern stabilisation* is achieved when $h >= h_2$.

**h-gaps**

It is more convenient now to identify related gaps at different levels by subtracting appropriate multiples of $a_k$, and so we introduce the concept of the *m-gap*:

$x > n(h_0-1)$ is an *m-gap* if $x_h$ is a gap at levels $h < m$, but is not a gap at level $m > h_0-1$, where:

$$x_h = x + (h - (h_0-1))a_k$$

We can now think of the range of values $n(h_0-1) < x < (h_0-1)a_k$ as divided into three groups:

Values with $(h_0-1)$-representations (and which are never gaps): we call these *pre-filled* values;

Values which are m-gaps for some $h_0-1 < m <= h_2$;

Values where $x_h$ is a gap at level h for all h: we call these *permanent* gaps.

An m-gap is a gap that is "filled in" at level m; once a gap is filled in, it remains so since if $x_h$ has an h-representation, then $x_{h+1}$ clearly has an (h+1)-representation using $a_k$ as the additional addend.

We define m(h) as the number of gaps at level h, and $d(h) = m(h-1) - m(h)$ for $h \geq h_0$ is then the number of gaps that are filled in at level h; that is, d(h) is the number of h-gaps.

*Lemma 1*

    $x < (h_0-1)a_k$ if x is an h-gap.

*Proof*

    By definition, $x_h$ has an h-representation and so must be less than or equal to $ha_k$, so:

$$x_h = x + (h - (h_0-1))a_k \leq ha_k$$

$$\Rightarrow \quad x \leq (h_0-1)a_k$$

    But $(h_0-1)a_k$ has an $(h_0-1)$-representation, and so is not an h-gap; so $x < (h_0-1)a_k$ as required.

*Lemma 2*

    If x is an h-gap with:

$$x_h = c_k a_k + c_{k-1} a_{k-1} + ... + c_1$$

    then $c_k = 0$.

*Proof*

$$x_h = x + (h - (h_0-1))a_k = c_k a_k + c_{k-1} a_{k-1} + ... + c_1 \qquad \text{with } c_k + ... + c_1 \leq h$$

$$\Rightarrow \quad x_{h-1} = x + ((h-1) - (h_0-1))a_k = (c_k-1)a_k + c_{k-1}a_{-1} + ... + c_1 \qquad \text{with } (c_k-1) + ... + c_1 \leq h-1$$

    So if $c_k > 0$, $x_{h-1}$ has an (h-1)-representation, and so x is not an h-gap.

*Theorem 5*

    Let x be an h-gap with:

$$x_h = c_{k-1}a_{k-1} + ... + c_1$$

    Then for any $c_i > 0$, $x + (a_k - a_i)$ is an (h-1)-gap.

*Proof*

    We write $x' = x + (a_k - a_i)$, and show first that $x'_{h-1}$ has an (h-1)-representation, and second that $x'_m$ has no m-representation for any $m < h-1$; hence x' is an (h-1)-gap.

    $x'_{h-1}$ has an (h-1)-representation:

$$x'_{h-1} = x' + ((h-1) - (h_0-1))a_k$$

$$= x + a_k - a_i + ((h-1) - (h_0-1))a_k$$

$$= x + (h - (h_0-1))a_k - a_i$$

$$= x_h - a_i$$

$$= c_{k-1}a_{k-1} + ... + (c_i-1)a_i + ... + c_1$$

which is an (h-1)-representation since $c_i > 0$.

$x'_m$ has no m-representation for any $m < h-1$:

Suppose $x'_m$ has an m-representation for some $m < h-1$:

$$x'_m = d_k a_k + ... + d_1 \qquad \text{with } d_k + ... + d_1 <= m$$

Then:

$$x_{h-1} = x + ((h-1) - (h_0-1))a_k$$

$$= x' - (a_k - a_i) + ((h-1) - (h_0-1))a_k$$

$$= x' + a_i + ((h-2) - (h_0-1))a_k$$

$$= x'_{h-2} + a_i$$

$$= x'_m + ((h-2) - m)a_k + a_i$$

$$= d'_k a_k + ... + d'_1 \qquad \text{where } d'_k = d_k + ((h-2) - m)$$
$$d'_i = d_i + 1$$
$$d'_j = d_j \text{ for values of j other than k, i}$$

So:

$$d'_k + ... + d'_1 = (d_k + ... + d_1) + ((h-2) - m) + 1 \ <= \ m + ((h-2) - m) + 1 = h-1$$

And so $x_{h-1}$ has an (h-1)-representation, contrary to our assumption that x is an h-gap.

**Selmer's erroneous proof that d(h) cannot increase as h increases**

Theorem 5 says that for each gap filled in at level h, there are one or more gaps filled in at level h-1; so the number of gaps filled in at level h can be no greater than the number filled in at level h-1.

*An example*

$A_4 = \{1, 4, 26, 35\}$; $h_0 = 8$, $h_1 = 9$, $h_2 = 10$, $n(h_0-1) = 22$, $(h_0-1)a_k = 245$

Here is a summary of the range of values $n(h_0-1) <= x <= (h_0-1)a_k$ where:

      `*` indicates a value that is pre-filled
      `-` indicates a value that is a permanent gap
      `d` indicates a value that is an h-gap: `8` for an 8-gap, `9` for a 9-gap, `0` for a 10-gap

```
        0         10        20        30        40        50        60        70
        0123456789012345678901234567890123456789012345678901234567890123456789
  0                       *8***************************************************
 80 *******************************************9*******8*********************8**-0****
160 ***-9********-8*******-***--*---***--*---***--*---***--*---*--------*--------*---
240 -----*
```

Note that we can determine both $h_1$ and $h_2$ from the list above:

$h_1$ is the smallest value of h such that the first gap is a permanent gap: that is, the gap at 154. Each increment in h from then on will add exactly $a_k$ to the h-range, because the cover will always be terminated by this gap.

We see that for $h = h_0-1 = 7$, there are four gaps below 154 that are yet to be filled in: 8-gaps 23, 129 and 151 and 9-gap 120. Three of these are filled in at h = 8, and the final one at h = 9: so $h_1 = 9$.

In general, $h_1$ is the largest value of h for which there is an h-gap below the first permanent gap.

$h_2$ is the smallest value of h for which all remaining gaps are permanent, and so is equal to the largest value of h for which there is an h-gap in the entire range; in this case we see that $h_2 = 10$.

Details of the h-gaps are as follows:

```
    h      x     x_h    h-representation

   10    155    260*   [0,10, 0, 0]**

    9    120    190    [0, 7, 2, 0]
    9    164    234    [0, 9, 0, 0]

    8     23     58    [0, 2, 0, 6], [0, 2, 1, 2]
    8    129    164    [0, 6, 2, 0]
    8    151    186    [0, 7, 1, 0]
    8    173    208    [0, 8, 0, 0]

    *    x_h = x + (h - (h_0-1))a_k = 155 + 3.35 = 260
    **   The representation  c_4a_4 + c_3a_3 + c_2a_2 + c_1 is denoted [c_4, c_3, c_2, c_1]
```

The relationships between the h-gaps are as follows (where the symbol -d> indicates that the h-gap on the right hand side is derived from that on the left by adding ($a_k - a_d$)):

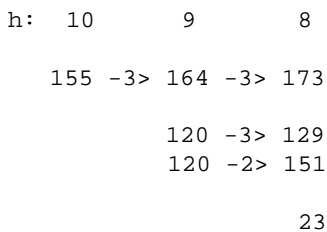

```
  h:   10        9         8

     155  -3>  164  -3>  173

               120  -3>  129
               120  -2>  151

                          23
```

Another way of illustrating the h-gaps and their relationships is by means of a diagram:

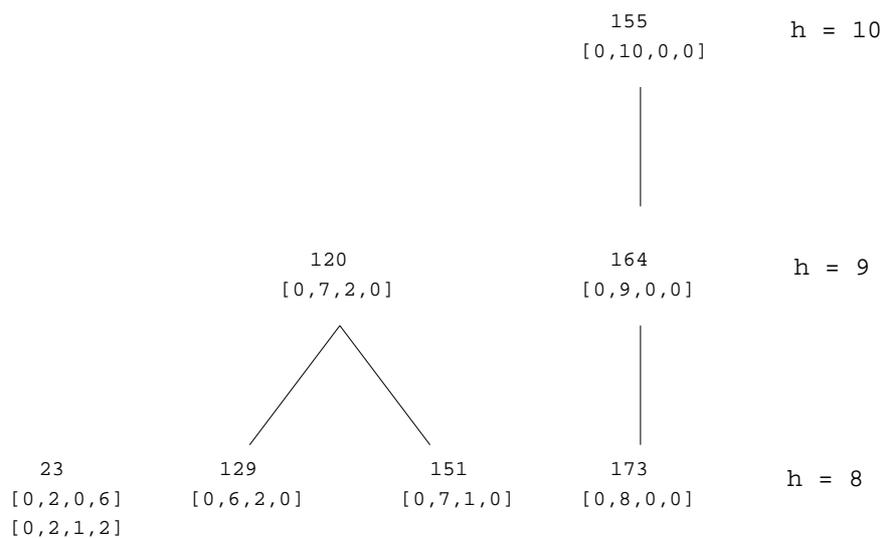

We see that d(8) = 4, d(9) = 2 and d(10) = 1 - all in accord with Selmer's conjecture.

Notes:

Since the representation* for the 10-gap 155 has only one non-zero coefficient, only one 9-gap can be derived; whereas the 9-gap 120 determines two 8-gaps since its representation has two non-zero coefficients.

The 8-gap 23 has two representations, one of which requires only 5 addends. This does not mean that it is an n-gap for some n < 8: there is no 7-representation for 23.

* For convenience, we talk of the "representation of the h-gap x" to mean "the h-representation of $x_h$"

*Another example*

$A_6 = \{1, 5, 8, 25, 37, 56\}$; $h_0 = 4$, $h_1 = 5$, $h_2 = 6$, $n(h_0-1) = 3$, $(h_0-1)a_k = 168$

We find:

```
    0         10        20        30        40        50        60        70
    0123456789012345678901234567890123456789012345678901234567890123456789

  0    *4********4******44*44****45**4***4***-***4***45**-*6****-4********4**-*5**---*
 80 4**---**-*--4**---**-*----*----***---**-*---------*------*-----------*----------
160 --------*
```

Details of the h-gaps are as follows:

```
    h    x    x_h   h-representation

    6    54   222   [0, 6, 0, 0, 0, 0]

    5    29   141   [0, 3, 1, 0, 1, 0]
    5    49   161   [0, 3, 2, 0, 0, 0]
    5    73   185   [0, 5, 0, 0, 0, 0]

    4     4    60   [0, 0, 2, 0, 2, 0]
    4    12    68   [0, 1, 1, 0, 1, 1]
    4    19    75   [0, 1, 1, 1, 1, 0], [0, 2, 0, 0, 0, 1], [0, 0, 3, 0, 0, 0]
    4    20    76   [0, 2, 0, 0, 0, 2], [0, 0, 3, 0, 0, 1]
    4    22    78   [0, 1, 1, 2, 0, 0]
    4    23    79   [0, 2, 0, 0, 1, 0]
    4    28    84   [0, 2, 0, 0, 2, 0]
    4    32    88   [0, 1, 2, 0, 0, 1]
    4    36    92   [0, 1, 2, 0, 1, 0]
    4    44   100   [0, 2, 1, 0, 0, 1], [0, 0, 4, 0, 0, 0]
    4    48   104   [0, 2, 1, 0, 1, 0]
    4    60   116   [0, 3, 0, 0, 1, 0]
    4    68   124   [0, 2, 2, 0, 0, 0]
    4    80   136   [0, 3, 1, 0, 0, 0]
    4    92   148   [0, 4, 0, 0, 0, 0]
```

We see that $d(4) = 15$, $d(5) = 3$, $d(6) = 1$ - again in accord with Selmer's conjecture.

The relationships between the h-gaps are as follows:

```
    h:  6          5          4

        54 -5>    73 -5>    92

                  29 -5>    48
                  29 -4>    60
                  29 -2>    80
                  49 -5>    68
                  49 -4>    80

                             4
                            12
                            19
                            20
```

```
                    22
                    23
                    28
                    32
                    36
                    44
```

with the corresponding diagram:

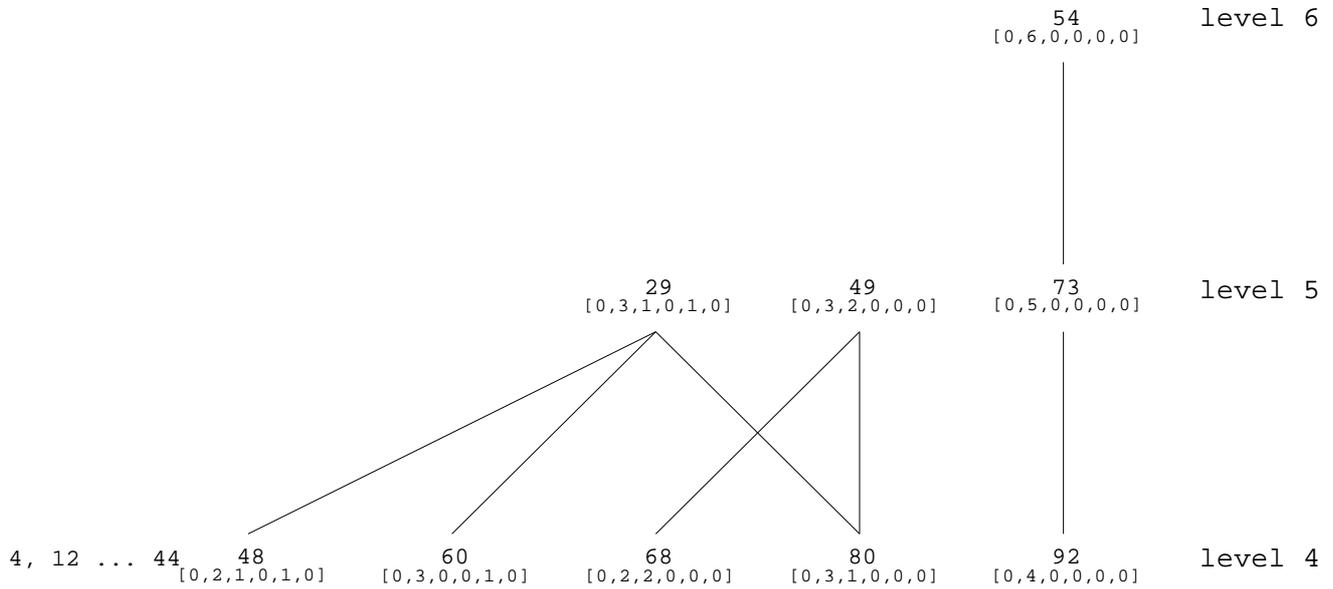

Notes:

> Here we see two different h-gaps determining the *same* (h-1)-gap: 5-gaps 29 and 49 both determine the 4-gap 80 (as well as other 4-gaps). This is why Selmer's argument is invalid: this possibility was not considered.

**My erroneous proof that d(h) cannot increase as h increases**

After failing to discover any counter-example to the conjecture despite extensive computing (see later), I became convinced that the conjecture was correct and started to search for a proof; this is what I came up with.

It seemed to me that the only way in which the number of h-gaps could decrease as h increased was if two different h-gaps determined the same (h-1)-gap *and no others*.

Now the only way in which one h-gap can determine only one (h-1)-gap is if its h-representation has only a single non-zero coefficient; for if it has at least two non-zero coefficients $c_i$ and $c_j$ then it must determine two (h-1)-gaps $x + (a_k - a_i)$ and $x + (a_k - a_j)$ (which cannot be the same since $a_i$ differs from $a_j$).

So we are looking for two h-gaps x and y with representations of that sort which both determine the same (h-1)-gap z *and no other*; we require:

$x_h = (c_p+1)a_p$ determines $z_{h-1} = c_p a_p$ where $z = x + (a_k - a_p)$
$y_h = (d_q+1)a_q$ determines $z_{h-1} = d_q a_q$ where $z = y + (a_k - a_q)$

Clearly p differs from q - for if p = q we find $c_p = d_q$ (from the formulae for $z_{h-1}$) and hence $x_h = y_h$.

Here's an (imaginary) illustration of the kind of relationship we are looking for:

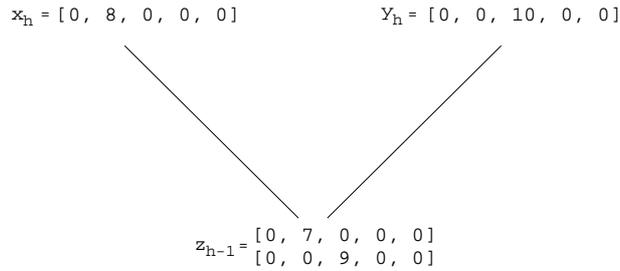

```
xh = [0, 8, 0, 0, 0]   -4>   zh-1 = [0, 7, 0, 0, 0]
yh = [0, 0,10, 0, 0]   -3>   zh-1 = [0, 0, 9, 0, 0]
```

Now $x_h = (c_p+1)a_p = c_p a_p + a_p = d_q a_q + a_p$  (via $z_{h-1}$)

Now $d_q + 1 <= h$, so this is another h-representation of $x_h$.

By Theorem 5, v is an (h-1)-gap where:

$v = x + (a_k - a_q)$       (note that $d_q > 0$, since $z_{h-1} > 0$)

Clearly v is different from z, since z - v = x - y; so our assumption that x determines only one (h-1)-gap is false.

In a similar way, we see that y must also determine a further (h-1)-gap, say $w = y + (a_k - a_p)$.

Our (still imaginary) illustration now looks like this:

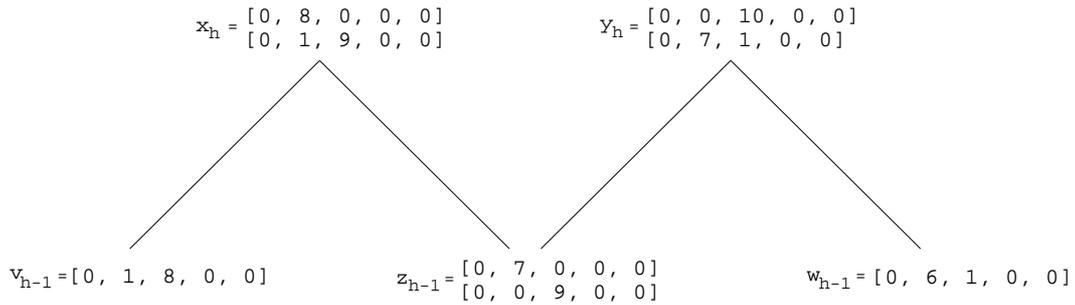

*Example showing equality of d(h) as h increases:*

$A_4 = \{1, 5, 18, 19\}$; $h_0 = 6$, $h_1 = 13$, $h_2 = 13$, $n(h_0-1) = 13$, $(h_0-1)a_k = 95$

```
     0         10        20        30        40        50        60        70
     01234567890123456789012345678901234567890123456789012345678901234567890123456789

 0              *6******************6*******************6******************9876********
80  **32109876******
```

We find $d(6, ... 13) = \{5, 2, 2, 2, 1, 1, 1, 1\}$:

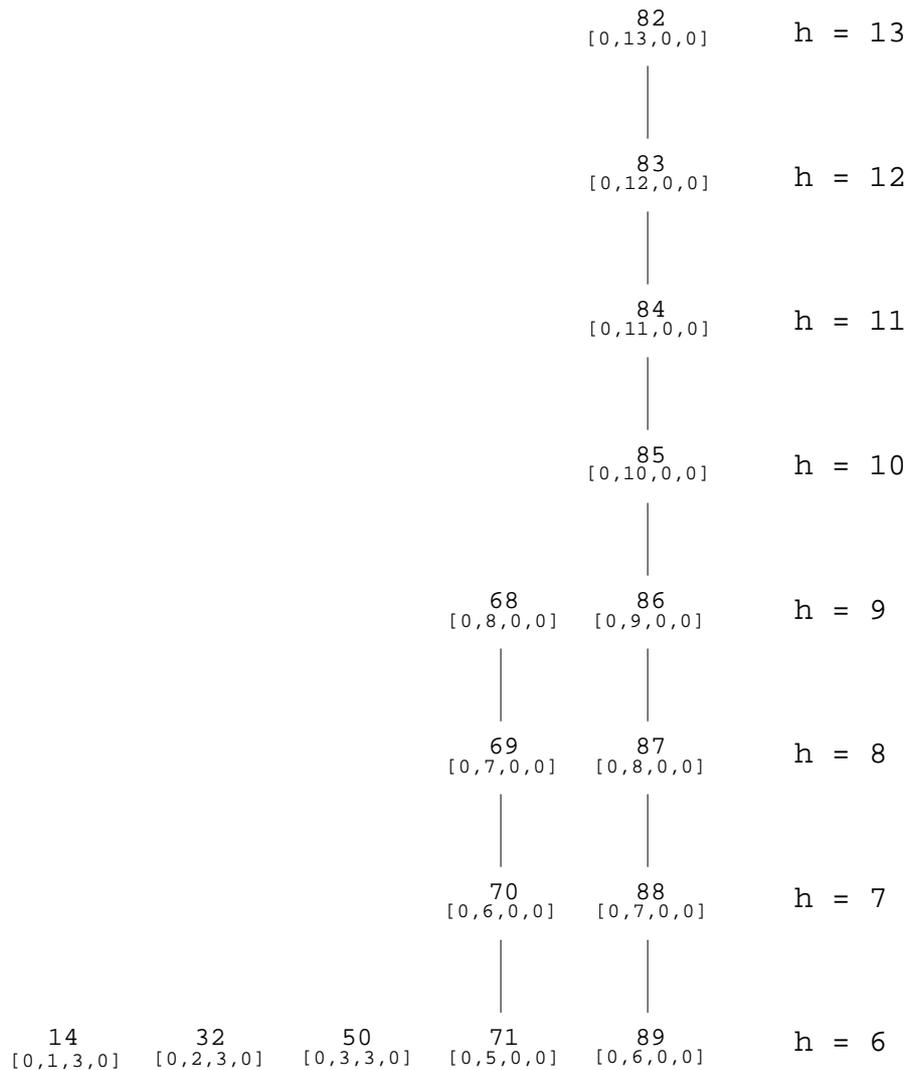

Note how every derived h-gap has a representation of the form [0, d, 0, 0], as expected; and that it is possible for two different h-gaps to have the same representation:

    8-gap 69 and 7-gap 88 both have the representation $7a_3 = 7.18 = 126$, since $69 + 3.19 = 88 + 2.19 = 126$.

*Example illustrating the case considered in my erroneous proof:*

$A_6 = \{1, 3, 8, 21, 28, 29\}$; $h_0 = 4$, $h_1 = 8$, $h_2 = 8$, $n(h_0-1) = 12$, $(h_0-1)a_k = 87$

```
         0         10        20        30        40        50        60        70
         01234567890123456789012345678901234567890123456789012345678901234567890123456789
   0              *4*4**4*4*****4****************4****654*****54******4****654**87654***
  80 7654****
```

We find $d(4, ... 8) = \{12, 5, 4, 2, 1\}$:

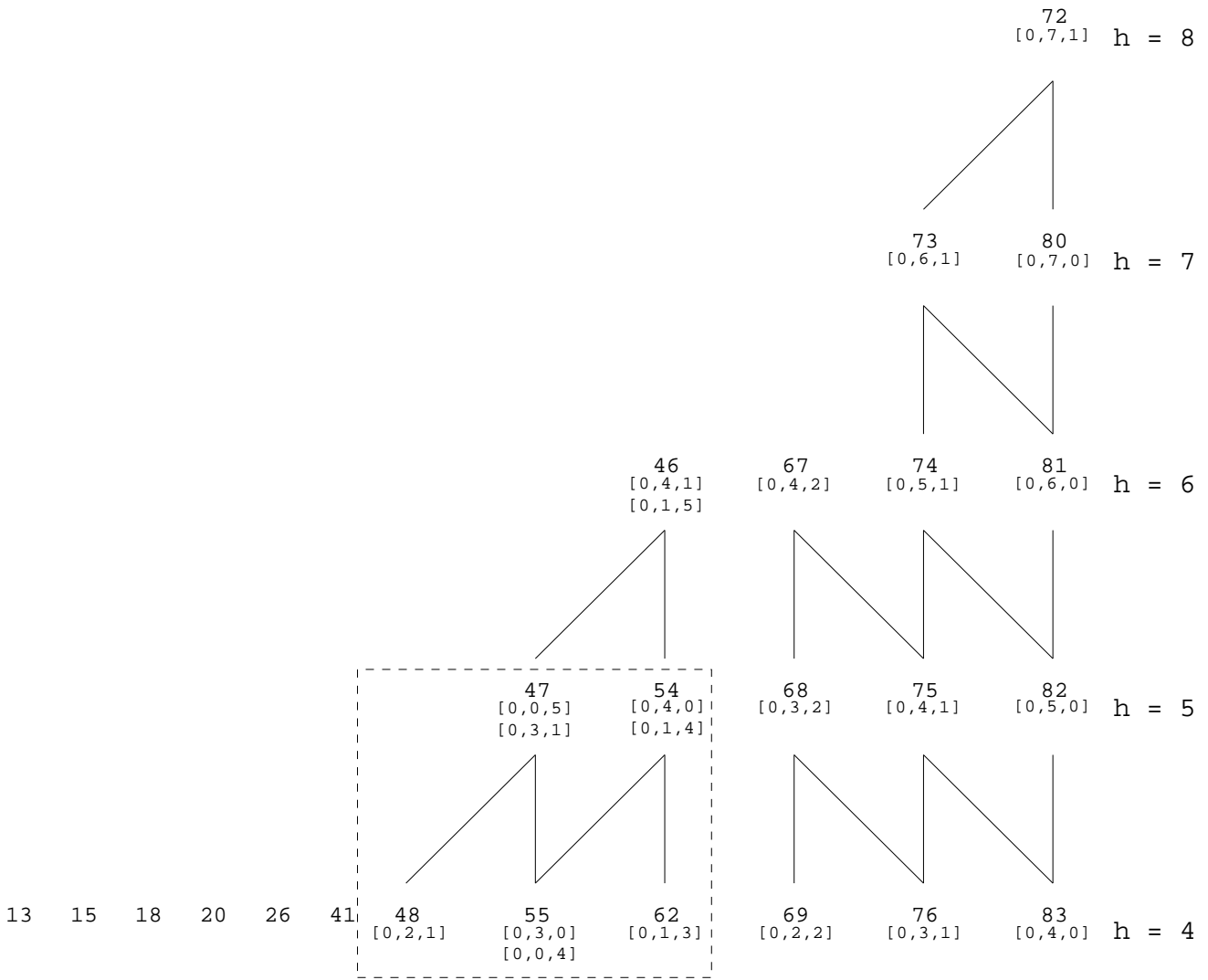

The structure inside the dotted rectangle is the one of interest (note that the representations are abbreviated in the diagram above: coefficients for $c_3$, $c_2$ and $c_1$ are always zero).

**A counter-example**

The problem with my "proof" is that - like Selmer - I do not consider all possible ways in which n h-gaps might determine (n-1) (h-1)-gaps. For example, there is no reason why three h-gaps might not determine two (h-1)-gaps as follows:

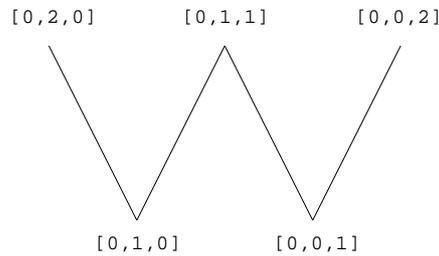

The first counter-example that I found is much more complex than this (see below), but this later example has many similarities:

$A_{10}$ = {1, 3, 4, 9, 12, 13, 19, 44, 47, 62}; $h_0$ = 3, $h_1$ = 5, $h_2$ = 5, $n(h_0-1)$ = 10, $(h_0-1)a_k$ = 124

```
            0         10        20        30        40        50        60        70
            0123456789012345678901234567890123456789012345678901234567890123456789

    0              *3**************3*33**35335*35335**5**5**4*-4**4**4**4**---3*-3**3--3
```

```
   80  -*------*--*--*------------*--*--------------*
```

We find d(3,4,5) = {14, 5, 6}: that is, d(5) > d(4), contrary to the conjecture:

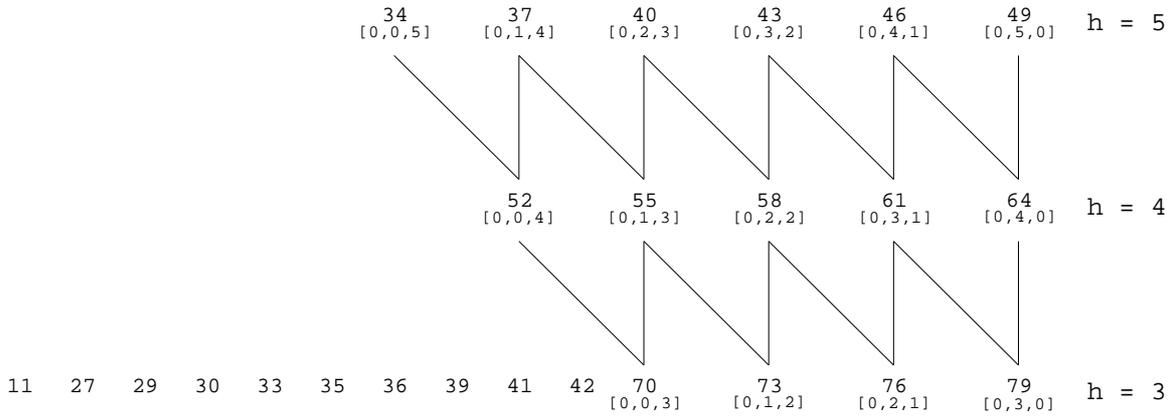

Note how similar is the pattern to the (imaginary) example above.

**Summary of results**

There are no exceptions to the conjecture for the following values of k and $h_0$:

| $h_0$ | k |
|---|---|
| 3 | 4, 5, 6, 7, 8 |
| 4 | 4, 5, 6, 7, 8 |
| 5 | 4, 5, 6, 7 |
| 6 | 4, 5, 6 |
| 7 | 4, 5, 6 |

| k | $h_0$ |
|---|---|
| 4 | <= 29 |
| 5 | <= 13 |

[k = 3 is not relevant, since then either $h_2 = h_0-1$ or $h_2 = h_0$ and there is at most one value d(h).]

There are seven exceptions for k = 9, $h_0$ = 3:

$A_9$ = {1, 2, 4, 6, 9, 10, 31, 32, 36}, $h_0$ = 3, $h_1$ = 7, $h_2$ = 7, n($h_0$-1) = 16, d(3 .. 7) = {9, 10, 6, 2, 1}
$A_9$ = {1, 2, 5, 8, 9, 10, 31, 33, 36}, $h_0$ = 3, $h_1$ = 6, $h_2$ = 6, n($h_0$-1) = 20, d(3 .. 6) = {9, 10, 5, 2}
$A_9$ = {1, 2, 5, 8, 9, 10, 31, 34, 36}, $h_0$ = 3, $h_1$ = 7, $h_2$ = 7, n($h_0$-1) = 20, d(3 .. 7) = {9, 10, 5, 3, 1}
$A_9$ = {1, 3, 4, 6, 9, 10, 31, 32, 36}, $h_0$ = 3, $h_1$ = 6, $h_2$ = 6, n($h_0$-1) = 16, d(3 .. 6) = {9, 10, 6, 2}
$A_9$ = {1, 3, 5, 6, 9, 10, 31, 32, 36}, $h_0$ = 3, $h_1$ = 6, $h_2$ = 6, n($h_0$-1) = 16, d(3 .. 6) = {9, 10, 6, 2}
$A_9$ = {1, 3, 5, 7, 9, 10, 31, 32, 36}, $h_0$ = 3, $h_1$ = 6, $h_2$ = 6, n($h_0$-1) = 20, d(3 .. 6) = {9, 10, 5, 1}
$A_9$ = {1, 3, 5, 7, 9, 10, 31, 35, 36}, $h_0$ = 3, $h_1$ = 8, $h_2$ = 8, n($h_0$-1) = 20, d(3 .. 8) = {9, 10, 6, 3, 2, 1}

In each case, the exception is the same: d(4) = 10 > d(3) = 9.

The first of these - discovered while looking for convincing illustrations of my "proof" - looks like this:

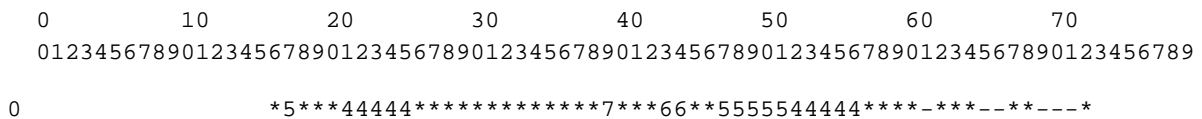

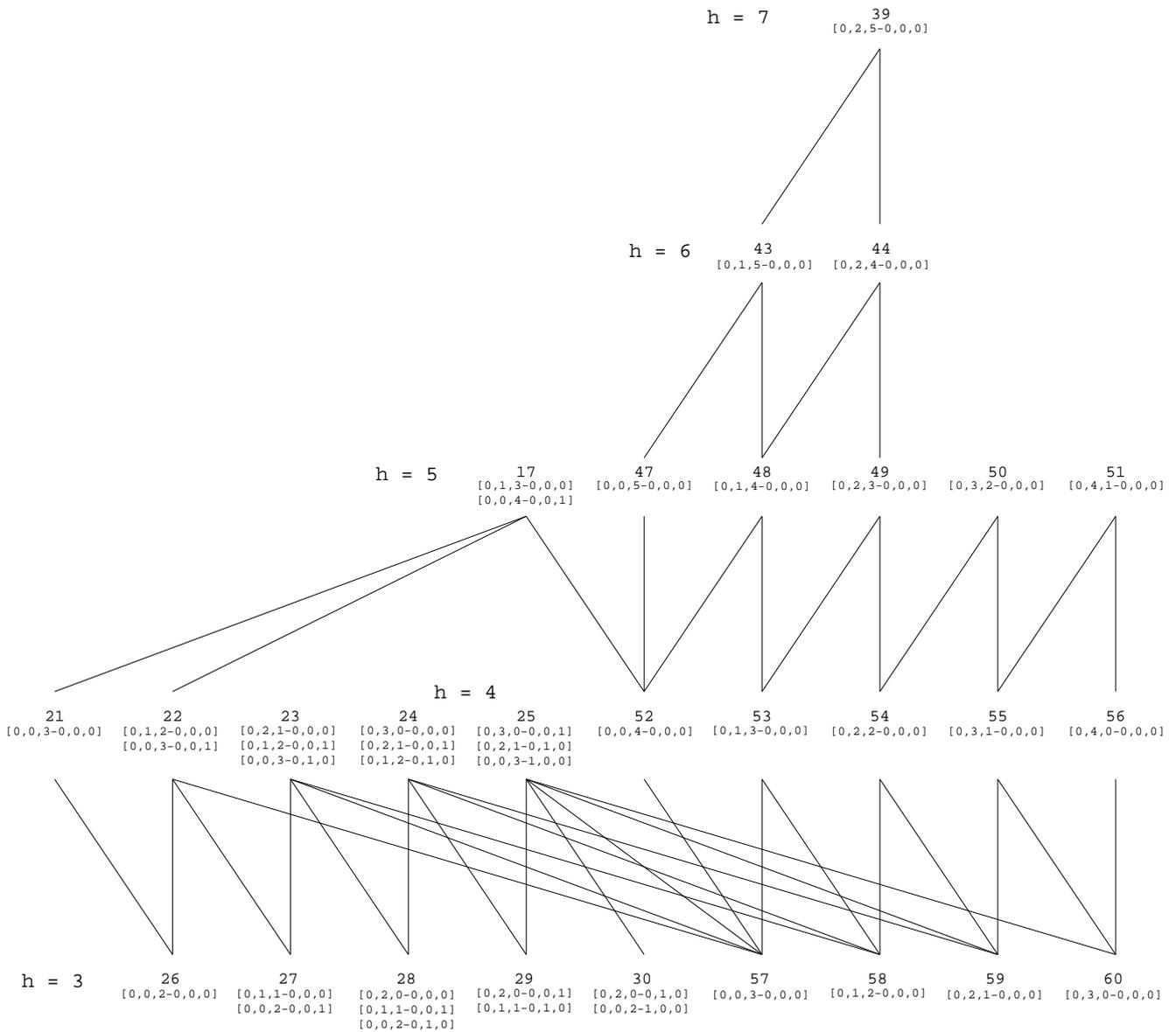

Note that each h-gap's representation is shown as [$c_9, c_8, c_7 - c_3, c_2, c_1$], since $c_6 = c_5 = c_4 = 0$.

For $k = 10$, $h_0 = 3$, there are 1031 exceptions of which all but four have $h_1 = h_2$; the odd ones out are:

$A_{10} = \{1, 3, 5, 12, 13, 15, 34, 53, 70, 75\}$, $h_0 = 3$, $h_1 = 8$, $h_2 = 11$, $n(h_0-1) = 6$, $d(3 .. 11) = \{26, 9, 10, 6, 4, 3, 1, 1, 1\}$
$A_{10} = \{1, 3, 7, 8, 12, 15, 18, 43, 54, 58\}$, $h_0 = 3$, $h_1 = 6$, $h_2 = 9$, $n(h_0-1) = 4$, $d(3 .. 9) = \{17, 8, 9, 4, 2, 1, 1\}$
$A_{10} = \{1, 3, 7, 8, 12, 18, 19, 43, 54, 58\}$, $h_0 = 3$, $h_1 = 6$, $h_2 = 9$, $n(h_0-1) = 4$, $d(3 .. 9) = \{20, 7, 8, 3, 1, 1, 1\}$
$A_{10} = \{1, 3, 7, 12, 15, 18, 19, 43, 54, 58\}$, $h_0 = 3$, $h_1 = 6$, $h_2 = 9$, $n(h_0-1) = 10$, $d(3 .. 9) = \{18, 7, 8, 3, 1, 1, 1\}$

There are also four cases where the difference $d(h) - d(h-1)$ is greater than one:

$A_{10} = \{1, 3, 4, 6, 7, 21, 35, 50, 64, 67\}$, $h_0 = 3$, $h_1 = 16$, $h_2 = 16$, $n(h_0-1) = 14$, $d(3 .. 16) = \{17,19,10,5,3,2,2,2,2,2,1,1,1,1\}$
$A_{10} = \{1, 3, 4, 9, 10, 24, 38, 53, 67, 73\}$, $h_0 = 3$, $h_1 = 10$, $h_2 = 10$, $n(h_0-1) = 14$, $d(3 .. 10) = \{18, 20, 8, 4, 3, 3, 2, 1\}$
$A_{10} = \{1, 3, 5, 7, 8, 24, 38, 57, 71, 73\}$, $h_0 = 3$, $h_1 = 10$, $h_2 = 10$, $n(h_0-1) = 16$, $d(3 .. 10) = \{20, 22, 13, 7, 6, 4, 2, 1\}$
$A_{10} = \{1, 4, 7, 8, 14, 15, 18, 48, 52, 55\}$, $h_0 = 3$, $h_1 = 9$, $h_2 = 9$, $n(h_0-1) = 2$, $d(3 .. 9) = \{15, 17, 11, 8, 4, 2, 1\}$

There is one case with the largest value of $h_2 = 18$:

$A_{10} = \{1,3,5,7,8,17,36,50,67,69\}$, $h_0 = 3$, $h_1 = 18$, $h_2 = 18$, $n(h_0-1) = 18$, $d(3 .. 18) = \{18,19,11,6,4,3,2,2,2,2,1,1,1,1,1,1\}$

There are twenty cases with the smallest value of $h_2 = 5$; here is one example:

$A_{10} = \{1, 3, 4, 9, 12, 13, 19, 44, 47, 62\}$, $h_0 = 3$, $h_1 = 5$, $h_2 = 5$, $n(h_0-1) = 10$, $d(3 .. 5) = \{14, 5, 6\}$

There are no cases where there is more than one occurrence of $d(h) > d(h-1)$ for different values of h.

**Miscellanea**

This is a list of the sets with highest $h_2$ value for various $k/h_0$:

| k | $h_0$ | $h_2$ | $h_1$ | $A_k$ |
|---|---|---|---|---|
| 4 | 10 | 35 | 35 | {1, 6, 41, 42} |
| 5 | 5 | 29 | 29 | {1, 6, 7, 36, 37} |
| 5 | 7 | 58 | 58 | {1, 8, 11, 69, 70} |
| 6 | 4 | 31 | 31 | {1, 5, 9, 10, 41, 42} |
| 7 | 3 | 23 | 23 | {1, 3, 8, 10, 11, 34, 35} |
|   |   |   |   | {1, 4, 7, 10, 11, 34, 35} |

This is a list of the sets with highest $(h_2-h_1)$ value for various $k/h_0$:

| k | $h_0$ | $h_2-h_1$ | $h_2$ | $h_1$ | $A_k$ |
|---|---|---|---|---|---|
| 4 | 10 | 4 | 15 | 11 | {1, 7, 39, 44} |
| 5 | 5 | 4 | 12 | 8 | {1, 6, 7, 34, 39} |
|   |   |   | 9 | 5 | {1, 6, 8, 35, 55} |
| 5 | 7 | 9 | 24 | 15 | {1, 7, 12, 64, 69} |
| 6 | 4 | 5 | 9 | 4 | {1, 2, 5, 13, 32, 43} |
|   |   |   | 11 | 6 | {1, 4, 9, 21, 37, 40} |
|   |   |   | 12 | 7 | {1, 4, 10, 11, 38, 43} |
| 7 | 3 | 5 | 8 | 3 | {1, 3, 5, 12, 13, 32, 29} |
|   |   |   | 7 | 2 | {1, 3, 7, 12, 22, 32, 43} |
| 7 | 4 | 10 | 16 | 6 | {1, 3, 10, 18, 22, 70, 73} |

**Programs**

The program called h0h1h2 scans all sets $A_k$ for a given value of $h_0$, looking for violations of the conjecture.

A secondary program, h012gaps, calculates $h_0$, $h_1$ and $h_2$ for a given set $A_k$, and then lists all values in the range $n(h_0-1, A_k) < x < (h_0-1)a_k$ indicating for each one whether it is an h-gap for some h, is a permanent gap, or is pre-filled.

**Finding all sets $A_k$ with a given $h_0$**

The algorithm used to iterate over all sets with a given value of $h_0$ took some time to get right, and is shown below. Although the version in h0h1h2 differs in some respects for efficiency reasons, it is easy to see that it is equivalent: the two functions have also been checked against each other experimentally.

```
void try (int i, int v)

 /*
    This function is called externally as try(2, 0) with a_1 = 1; the eventual
     effect is to call process(x) for each set A_k which is h-admissible only
     for h >= h_0, with x = n(h_0-1, A_k).

    On entry, a_1, a_2 ... a_{i-1} have been chosen, and A_{i-1} is known to be an
     h_0-admissible set.

    If all values x < a_{i-1} can be generated using less than h_0 stamps, v = 0
```

```
      (and the set is h-admissible for some value(s) h < h_0).

   Otherwise, v > 0 is the smallest value v < a_{i-1} which requires exactly h_0
    stamps to generate (and the set is h-admissible only for h >= h_0, with
    n(h_0-1, A_{i-1}) = v-1).

   Values a_i = a_{i-1}+1, a_{i-1}+2 ... max are considered in turn, where max is the
    smallest value that does not have an h_0-generation.

   If i < k, the recursive call try(i+1, v') is made for each value (where v'
    is the analogue of v for i+1).

   If i = k and v' > 0, the set is h-admissible only for h >= h_0 with
    v' = n(h_0-1, A_k) + 1, and process(v'-1) is called as required.
*/

{
   a[i] = a[i - 1] + 1;

   /* The body of this loop considers the set 1, a_2, ... , a_{i-1}, a_i.
      It is known that all values x <= a_i-2 are h_0-admissible. */

   do
   {
      int h = find_minimal (a[i] - 1, i - 1, h0);
       /* This call looks for the minimal h-generation for h <= h_0: if a_i-1 has any
            generation using at most h_0 stamps, then h will be the smallest number
            of stamps required to generate a_i-1; otherwise h = 0 */

      if (h == 0)     /* first gap; so no set with a_i greater than or equal to the */
          return;     /*  current value can ever be h_0-admissible */

      if (h == h0)    /* a critical value that requires exactly h_0 stamps */
      {
         /* if this is the first such value encountered, it defines n(h_0-1, A_k) to
            be equal to a_i-2 */
         if (v == 0)
             v = a[i] - 1;
      }

      /* A_i is h-admissible for some minimal h <= h_0, and, if this minimal value
         h = h_0, then v = n(h_0-1, A_i) + 1.

         If h = h_0, then the loop below simply processes the single value a_i.

         But if h < h_0, there may be further values a_i+1 ... a_i+(h_0-h-1) that we
          can process without further ado - since we know that if a_i-1 has an
          h-generation for h < h_0-1, then a_i has an h-generation for h < h_0. */
      do
      {
         if (i == k)
         {
            /* only call process(..) if A_k is h-admissible only for h >= h_0 */
            if (v > 0)
                process (v - 1);
         }
         else
            try (i + 1, v);

         a[i]++;
         h++;

      } while (h < h0);

   } while (TRUE);
}
```

This algorithm makes use of a function `find_minimal(v,i,s)` which in turn calls a function `can_minimal(v,i,s)` to find the h-generation (if any) of v which has minimal h for h <= s. This function differs from the standard function `can(v,i,s)` only in that it has to range over *all* possible generations instead of returning as soon as the first one has been found.

**Evaluating h-gap information - algorithm A**

We first find a value $H > h_2$ as follows.

By definition, no gap $n(h) < x_h < ha_k$ can ever be filled in for $h >= h_2$; that is, there is no h-representation

$$x_h = c_{k-1}a_{k-1} + ... + c_1 \text{ with } c_{k-1} + ... c_1 <= h \text{ for any } h >= h_2$$

This is certainly true as soon as $x_h > ha_{k-1}$: for each increment in h increases $x_h$ by $a_k$ and provides at best only one more addend $a_{k-1}$.

So we need to find H such that $x_H > Ha_{k-1}$; now $x_H > n(H)$, so a value H satisfying $n(H) >= Ha_{k-1}$ will do.

But by Theorem 1 we know $n(H) >= n(h_0-1) + (H - (h_0-1))a_k$, so the following is sufficient:

$$n(h_0-1) + (H - (h_0-1))a_k >= Ha_{k-1}$$

$$\Rightarrow \quad H(a_k - a_{k-1}) >= (h_0-1)a_k - n(h_0-1)$$

So the following value will do:

$$H = \text{intpt}[ ((h_0-1)a_k - n(h_0-1) - 1) / (a_k - a_{k-1}) ] + 1$$

Note:
 This derivation was made before re-reading Selmer's proof of the existence of $h_1$ and $h_2$ where the value $h_s$ is introduced.
 Using his function $D(h) = n(h) - ((h+1)a_{k-1} - a_k)$, a slightly better value of H is calculated as follows:
  $D(h) >= 0$
  $\Rightarrow \quad n(h) >= ((h+1)a_{k-1} - a_k)$
 which is true if: $n(h_0-1) + (h - (h_0-1))a_k >= ((h+1)a_{k-1} - a_k)$     (by Theorem 1)
  $\Rightarrow \quad h(a_k - a_{k-1}) >= (h_0-1)a_k + a_{k-1} - a_k - n(h_0-1)$
  $\Rightarrow \quad h(a_k - a_{k-1}) >= (h_0-2)a_k + a_{k-1} - n(h_0-1)$.

We can now proceed as follows to classify each value $n(h_0-1) < x < (h_0-1)a_k$:

If $x_H$ has no H-generation, then x is a permanent gap (since $H >= h_2$).

Otherwise, let:

$$x_H = x + (H - (h_0-1))a_k = c_k a_k + c_{k-1}a_{k-1} + ... + c_1 \quad \text{with } c_k + c_{k-1} + ... + c_1 <= H$$

be the H-generation with maximal coefficient* $c_k$, and let $m = H - c_k$ so that:

$$x_m = c_{k-1}a_{k-1} + ... + c_1 \quad \text{with } c_{k-1} + ... + c_1 <= m$$

* The standard can(v,i,s) algorithm naturally works through possibilities starting with the largest $c_k$ first, so only a simple modification to record the coefficient $c_k$ is necessary.

Note that $x_m$ cannot have another m-representation with $c_k > 0$, since that would contradict our assumption that the H-representation for $x_H$ has maximal coefficient $c_k$.

If $m <= h_0-1$, then $x_h$ has an h-generation for all $h >= h_0-1$ and so is a pre-filled value.

If $m > h_0-1$, then x has no h-generation for $h < m$ and thus is an m-gap.

To calculate $h_1$ and $h_2$ we proceed as follows:

We keep track of the maximum value of m for all m-gaps encountered as x increases.

$h_1$ is this maximum value as soon as the first permanent gap is encountered*.

$h_2$ is this maximum value over the entire range.

* If no permanent gap is encountered, $h_1 = h_2$

To calculate d(h) for $h_0 <= h <= h_2$ we proceed as follows:

An array delta[$h_0$ ... H] is allocated and initialised to zero.

Each time an m-gap is found, delta[m] is incremented.

When the whole range has been processed, d(h) is given by the value of delta[h].

Selmer's conjecture can now be tested directly by comparing the values in the array delta: a counter-example has been found if we discover delta[h] > delta[h-1] for some value of h.

The detailed workings of this algorithm were checked thoroughly by including independent checks on the various values calculated at appropriate points in the program:

Code to verify the value of $n(h_0-1)$ by direct calculation, using the standard can(x, k, $h_0$-1) function starting at x = 0 to find the first value that has no ($h_0$-1)-representation.

Code to verify the value of $h_0$ by direct calculation, continuing from x = $n(h_0-1)$ with can(x, k, $h_0$) to check that $n(h_0) > a_k$.

Code to verify the delta[h] vector: the number of gaps for each value of h is calculated by checking each value $x_h$ to see if it has an h-representation.

**Evaluating h-gap information - algorithm B**

It turns out that H can become surprisingly large, in which case the standard call can($x_H$, k, H) may take quite some time to determine that no H-representation exists (in the case of permanent gaps). In my first naive implementation a large H also led to a large vector delta[..] and correspondingly lengthy initialisation to zero; but this is easy to remedy: simply allocate a fixed size delta[..] large enough for all expected $h_2$, and initialise values only as they are needed.

Nonetheless, I sought an alternative algorithm which might be computationally more efficient, and came up with a new scheme based on Theorem 5, which states:

If x is an h-gap, then $x + (a_k - a_i)$ is an (h-1)-gap.

The algorithm processes values $n(h_0-1) < x < (h_0-1)a_k$ in sequence as follows, recording results in a pre-allocated array called gaps[x]:

Look for the $h_0$-representation of $x_{h0} = x + a_k$ which has maximum coefficient $c_k$.

If no such representation exists, mark x as a gap by setting gaps[x] to some large positive value (guaranteed to be larger than $h_2$); at this point we know only that x is not pre-filled - it might be an h-gap for some $h_2 >= h > h_0$, or it might be a permanent gap.

If a representation is found with $c_k > 0$, then mark x, x+1, ... x+nsp as pre-filled by setting gaps[y] = $h_0$-1 for x <= y <= x+nsp, where nsp is the number of spare stamps in the $h_0$-generation: nsp = $h_0 - (c_k + c_{k-1} + ... + c_1)$.

If a representation is found with $c_k = 0$, then x is an $h_0$-gap, and gaps[x] is accordingly set to $h_0$; the following algorithm is then applied recursively by calling propagate(x, $h_0$):

```
void propagate (int x, int hx)
{
```

```
            int i = k - 1;

            hx++;

            do
            {
                int y = x - (a_k - a_i);

                i--;

                if (y <= n(h_0-1))
                    return;

                if (gaps[y] > hx)
                {
                    gaps[y] = hx;
                    propagate (y, hx);
                }

            } while (i > 0);
        }
```

In effect, the function propagate (..) implements the inverse of Theorem 5, and when the entire range has been processed in this way we have:

$\text{gaps}[x] = h_0-1$      if x is a pre-filled value;

$\text{gaps}[x] = \text{large number}$      if x is a permanent gap;

$\text{gaps}[x] = h$      if x is an h-gap.

We show that this is the case as follows:

First, consider the effects of the recursive function propagate(x, $h_x$).

We claim that at the time of the call, x is known to be an h-gap for some h <= hx. The algorithm then looks at each value $y = x - (a_k - a_i)$, $y > n(h_0-1)$ in turn. We have:

$$x_{hx} = x + (hx - (h_0-1))a_k \quad \text{with hx-representation } c_k a_k + ... + c_i a_i + ... + c_1$$

and so:

$$\begin{aligned} y_{hx+1} &= y + (hx+1 - (h_0-1))a_k \\ &= x - (a_k - a_i) + (hx+1 - (h_0-1))a_k \\ &= x + (hx - (h_0-1))a_k + a_i \\ &= x_{hx} + a_i \end{aligned}$$

has (hx+1)-representation $c_k a_k + ... + (c_i+1)a_i + \ + c_1$, and so y is an h-gap for some h <= (hx+1).

But the initial call propagate(x, $h_0$) satisifies the condition that x is an h-gap for some h <= $h_0$ (actually equality holds), and so our claim is proved by induction.

For each such y we look at gaps[y].

If gaps[y] = a large number, then we set gaps[y] = hx+1 to indicate that we now know that y is not a permanent gap: instead it is an h-gap for some h <= hx+1. We then recurse.

If gaps[y] = $h_0-1$, we know that y is a pre-filled value, and so there is no need to recurse further.

If gaps[y] > hx+1, we set it to hx+1 - thus recording the new upper limit on the gaps's order - and recurse.

If $h_0$-1 < gaps[y] < hx+1 we already have a better bound for the order of the gap, and must have recursed over this value previously, so there is no need to process further.

So when all values $n(h_0-1) < x < (h_0-1)a_k$ have been processed, we know that:

gaps[x] = $h_0$-1             if x is a pre-filled value

gaps[x] = hx > $h_0$-1        if x is either a permanent gap or an h-gap for some h <= hx

We now look at it from the other side, and show that if x is an h-gap, then gaps[x] = h.

By Theorem 5, x is an h-gap => $x + (a_k - a_i)$ is an (h-1)-gap for at least one value of i (where $c_i > 0$); applying this repeatedly we see that there must be at least one value y derived from x via (h - $h_0$) steps that is an $h_0$-gap. But we have seen that propagate(y, $h_0$) will mirror these steps in reverse order, and so will set gaps[x] = h as required.

Comparison with Algorithm A:

Each value x is processed by a call can($x+a_k$, k, $h_0$) rather than can($x_H$, k, H) - which should mean a substantial time saving in many cases. On the other hand, each $h_0$-gap has to be propagated to determine where all of the h-gaps reside.

Algorithm B requires extra storage: the gaps[..] vector is approximately $h_0 a_k$ in length.

In practice, it turns out that there is not much to choose between the two algorithms, although there is some evidence to suggest that B is better than A for small $h_0$, large k, and A is better than B for small k, large $h_0$:

| k | $h_0$ | $t_A$ | $t_B$ | B v A |
|---|---|---|---|---|
| 4 | 12 | 15.68 | 15.80 | -0.8% |
| 4 | 15 | 83.27 | 85.03 | -2.1% |
| 5 | 7 | 85.51 | 83.54 | +2.3% |
| 7 | 3 | 14.34 | 13.66 | +4.7% |
| 8 | 3 | 389.45 | 366.03 | +6.0% |

(Times in seconds on a RISC PC with StrongARM processor)

Other interesting statistics for algorithm B:

| k | h | (a) | (b) |
|---|---|---|---|
| 4 | 12 | 1.11 | - |
| 4 | 15 | 1.46 | - |
| 5 | 7 | 2.30 | 0.03% |
| 7 | 3 | 1.93 | 1.34% |
| 8 | 3 | 3.03 | 1.98% |

(a) is the average number of elements of the gaps vector filled in by propagate (..) per set considered. I find this number surprisingly small: one possible explanation is that the real h-gaps are identified early on, and so most calls to propagate (..) do not recurse at all.

(b) is A/B as a percentage, where B is the number of times propagate (..) fills in an element of the gaps vector, and A is the number of times that that is an "improvement" over the previous value. For example, suppose gaps[x] = 10 (that is, x is known to be an h-gap for some h <= 10) and propagate (..) decides to set gaps[x] = 9 - that's an "improvement". We see that this happens very rarely indeed - another reason why (a) is so low.

**Development history**

Once algorithm A had been tested and checked, I started the search for a counter-example to Selmer's conjecture. After a total of some 40 hours' computing on my RISC PC with StrongARM processor I became convinced that no such counter-examples existed, and that the conjecture was correct. By this time, I had covered the following ranges of $k/h_0$:

    4/3 to 4/29
    5/3 to 5/13
    6/3 to 6/7
    7/3 to 7/5

I then turned my attention to algorithm B for some time, and later made some modifications to it to record the relationships between h-gaps arising from Theorem 5. With this additional information represented in the computer, I was able to see if I could find any "local" violation of the conjecture. I started by looking at cases where one (h-1)-gap x was derived from two or more h-gaps $y_1, ... y_n$, and then counting how many different (h-1)-gaps were determined by the h-gaps $y_1$ to $y_n$; but I found no cases where there was a smaller number of (h-1)-gaps than h-gaps over the following ranges:

    4/3 to 4/22
    5/3 to 5/12
    6/3 to 6/6

Finally convinced that the conjecture must be true, I produced a "proof" instead, and then started searching for some real examples to illustrate some of the proof's aspects. Imagine my surprise when the program claimed to have found a counter-example: indeed, I did not at first believe it, and it was only after a thorough check by hand that I had to agree that my proof must be wrong and that the conjecture had at last been shown to be false!

Further runs showed a total of 3 counter-examples for 9/3, and numerous further ones for 10/3; I noticed that all of these had $h_1 = h_2$, and wondered if this would always be the case.

But the development of a false proof was not to be my final embarrassment with this project: while checking the results just before writing them up, I found to my horror that a simple but major error had been present in the try (..) function throughout my experiments. The effect of this was that no sets of the form $\{1, h_0, a_3, ... a_k\}$ had ever been processed! I hastily corrected the error and re-ran previous results to confirm that no counter-examples existed in the ranges covered above. I also re-ran 9/3 and 10/3, finding examples where $h_1$ differs from $h_2$ for the first time.